\documentclass{article}
     
\usepackage{graphicx}       
\usepackage{cite}      	                                 

\usepackage{amsmath}
\usepackage{amsfonts}
\usepackage{amssymb}
\usepackage{xcolor}
\usepackage{hyperref}

	\newcommand\R{\mathbb{R}}

	\newcommand\hNir{\mathcal{H}}
	\newcommand\hVis{\mathcal{H}_{\rm Vis}}
	\newcommand\hAlt{\mathcal{H}_{\rm Alt}}

	\newtheorem{thm}{Theorem}[section]
	
	\newtheorem{dfn}[thm]{Definition}

	\newtheorem{con}[thm]{Condition}

	\usepackage{subfig}
	\usepackage{caption}
	\usepackage{url}
	\usepackage{float}

    \makeatletter
	    \let\@fnsymbol\@arabic
    	\newcommand{\specificthanks}[1]{\@fnsymbol{#1}}
    \makeatother	
	
	\title{Recent Advances in Reaction-Diffusion Equations 	with Non-Ideal Relays}

	\date{}

	\author{Mark Curran\thanks{Free University of Berlin, Institute of Mathematics I, Arnimallee 7, 14195, Berlin Germany 	\texttt{mark.curran88@gmail.com}} \and Pavel Gurevich\thanks{Free University of Berlin, Institute of Mathematics I, Arnimallee 3, 14195, Berlin Germany \texttt{gurevich@math.fu-berlin.de}} \textsuperscript{,}\thanks{Peoples' Friendship University of Russia 117198, Moscow Miklukho-Maklaya str.~6} \and Sergey Tikhomirov\thanks{Max Planck Institute for Mathematics in the Sciences, Inselstra\ss e 22, 04103 Leipzig, Germany \texttt{sergey.tikhomirov@gmail.com}} \textsuperscript{,}\thanks{Chebyshev Laboratory, St. Petersburg State University, 14th Line, 29b, Saint Petersburg, 199178 Russia}}

\begin{document}

	\maketitle

	\begin{abstract}
			We survey recent results on reaction-diffusion equations with discontinuous
hysteretic nonlinearities. We connect these equations with free
boundary problems and introduce a related notion of spatial
transversality for initial data and solutions. We assert that the equation
with transverse initial data possesses a unique solution, which
remains transverse for some time, and also describe its regularity. At
a moment when the solution becomes nontransverse, we discretize the
spatial variable and analyze the resulting lattice dynamical system
with hysteresis. In particular, we discuss a new pattern formation
mechanism --- {\it rattling}, which indicates how one should reset the
continuous model to make it well posed.
	\end{abstract}

	\section{Introduction}
	\label{sec:intro}
	\subsection{Motivation}
		
		In this chapter we will survey recent results on \index{reaction-diffusion equations} reaction-diffusion equations with a \index{hysteretic discontinuity} hysteretic discontinuity defined at every spatial point.
		We also refer to \cite{Visintin-1994-differential,krejc-1996-hysteresis,BroSpr-1996-hysteresis} and the more recent surveys by Visintin \cite{Visintin-2014-ten,Visintin-2015-P.D.E.s} for other types of partial differential equations with hysteresis.
		
		The equations we are dealing with in the present chapter were introduced in \cite{HopJae-1980-pattern, HopJaePoe-1984-hysteresis} to describe growth patterns in colonies of bacteria (Salmonalla typhirmurium). In these experiments, bacteria (non-difussing) are fixed to the surface of a petri dish, and their growth rate responds to changes in the relative concentrations of available nutrient and a growth-inhibiting by-product. The model asserts that at a location where there is a sufficiently high amount of nutrient relative to by-product, the bacteria will grow. This growth will continue until the production of by-product and diffusion of the nutrient lowers this ratio below a lower threshold, causing growth to stop. Growth will not resume until the diffusion of by-product raises the relative concentrations above an upper threshold that is distinct from the lower. Numerics in \cite{HopJae-1980-pattern} reproduced the formation of distinctive concentric rings observed in experiments, however the question of the existence and uniqueness of solutions, \index{existence and uniqueness of solutions} as well as a thorough explanation of the mechanism of pattern formation, \index{pattern formation} remained open.
		
		Another application in developmental biology can be found, e.g., in \cite{Marciniak-2006-receptor}, and an analysis of the corresponding stationary solutions in \cite{Koethe-2013-hysteresis}.

	\subsection{Setting of the Problem}
	
\begin{figure}
      \includegraphics[width=1.0\linewidth]{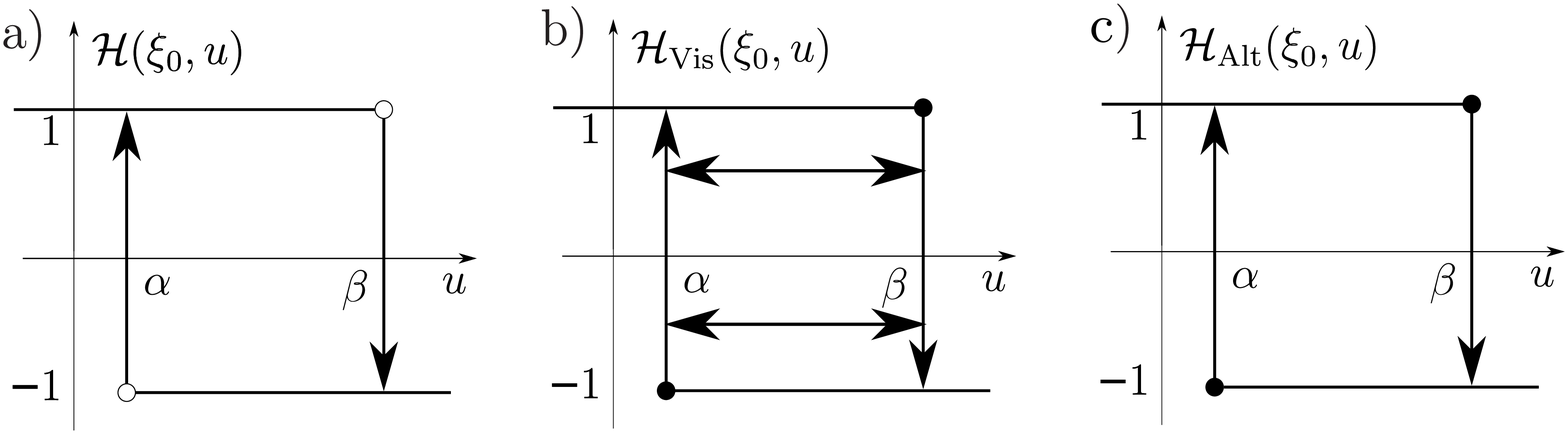}
			\caption{The hysteresis operator with $H_1(u) \equiv 1$ and $H_{-1} (u) \equiv -1$.}\label{fig:hyst-output}
			\end{figure}

		In this chapter we will treat the following prototype problem:
		\begin{align}
			& u_t = \Delta u + f(u,v), \quad v = \hNir(\xi_0,u), \quad (x,t) \in Q_T,  \label{problem-eq1} \\
			& u|_{t=0} = \varphi, \quad x \in Q, \label{problem-eq2} \\
			& \frac{\partial u}{\partial \nu} \bigg|_{\partial^\prime Q_T} = 0. \label{problem-eq3}
		\end{align}
		Here $Q \subset \R^n$ is a domain with smooth boundary, $Q_T := Q \times (0,T)$, where $T>0$, $\partial^\prime Q_T := \partial Q \times (0,T)$, $u$ is a real-valued function on $Q_T$, and $\mathcal{H}(\xi_0,u)$ is a hysteresis operator defined as follows (see Figure \ref{fig:hyst-output}a).
		Fix two real numbers ${\alpha < \beta}$, an integer $\xi_0 \in \{ -1, 1 \}$, and two continuous functions $H_1: (-\infty,\beta] \rightarrow \R$ and $H_{-1} : [\alpha,\infty) \rightarrow \R$ such that $H_1(u)\ne H_{-1}(u)$ for $u\in[\alpha,\beta]$.
		Define the sets
		\[ \Sigma_{1} := \{ (u,v) \in \R^2 \mid u \in ( -\infty,\beta), v= H_1(u) \}, \]
		\[ \Sigma_{-1} := \{ (u,v) \in \R^2 \mid u \in ( \alpha,\infty), v= H_{-1}(u) \}. \]
			\begin{dfn}\label{dfn-hNir}
				Let $u,v: [0,T] \rightarrow \R$, where $u$ is a continuous function. We say that $v = \hNir (\xi_0,u)$ if the following hold{\rm :}
				\begin{enumerate}
					\item $(u(t),v(t)) \in \Sigma_{1} \cup \Sigma_{-1}$ for every $t \in [0,T]$.
					\item If $u(0) \in (\alpha,\beta)$, then $v(0) = H_{\xi_0} (u(0))$.
					\item\label{dfn-hNir-item3} If $u(t_0) \in (\alpha,\beta)$, then $v(t)$ is continuous in a neighorhood of $t_0$.
				\end{enumerate}	
			\end{dfn}
		The operator $\hNir(\xi_0,u)$ is called the \emph{non-ideal relay} \index{non-ideal relay} and item \ref{dfn-hNir-item3} means that the non-ideal relay jumps up (or down) when $u = \alpha$ (or $u = \beta$). This definition is equivalent to the definitions of non-ideal relay found in \cite{Krasnosel-2012-systems, Visintin-1994-differential,GurTikSha-2013-reaction}.
		If ${{\hNir(\xi_0,u) (t) = H_j (u(t))}}$, then we call $\xi(t) := j$ the configuration of $\hNir$ at the moment $t$, and we call $\xi_0$ the initial configuration.
		Now let $u : Q_T \rightarrow \R$ be a function of $(x,t)$ and $\xi_0: Q \rightarrow \{-1,1\}$ a function of $x$, then $\hNir(\xi_0,u)(x,t)$ is defined in the same way by treating $x$ as a parameter, i.e., there is a non-ideal relay at every $x \in Q$ with input $u(x,t)$, configuration $\xi(x,t)$, and initial configuration $\xi_0(x)$.
	\subsection{Set-Valued Hysteresis}
		First results on the well-posedness\index{well-posedness} of \eqref{problem-eq1}--\eqref{problem-eq3} were obtained in \cite{Alt-1985-thermostat,Visintin-1986-evolution} for\index{set-valued hysteresis} set-valued hysteresis, and their model problems are worth explaining in more detail. In both papers, the \index{uniqueness of solutions} uniqueness of solutions as well as their continuous dependence on initial data remained open.
				
		First we discuss the work of Visintin \cite{Visintin-1986-evolution}, which treats \eqref{problem-eq1}--\eqref{problem-eq3} for arbitrary $n \geq 1$ with $\mathcal{H}(\xi_0,u)$ replaced by a set-valued operator called a \emph{completed relay}\index{completed relay} (see Figure \ref{fig:hyst-output}b).
		We still use the thresholds $\alpha < \beta$, and will consider constant hysteresis branches $H_1 (u) \equiv 1$, and $H_{-1} (u) \equiv -1$.
		We also define the set $\Sigma_{0} := \{ (u,v) \in \R^2 \mid u \in [\alpha,\beta], v \in (-1,1) \}$.
		
		\begin{dfn}\label{dfn-kVis}
			Let $u,v : [0,T] \rightarrow \R$, where $u$ is a continuous function, and let $\xi_0 \in [-1,1]$.			
			We say $v \in \hVis (\xi_0,u)$ if the following hold{\rm :}
				\begin{enumerate}
				
					\item $(u(t),v(t)) \in \overline{\Sigma_1} \cup \overline{\Sigma_{-1}} \cup \Sigma_0$ for every $t \in [0,T]$.

					\item\label{dfn-Kvist-item2} If $u(0) \in (\alpha, \beta)$, then $v(0) = \xi_0${\rm ;} if $u(0) = \alpha$ {\rm(}or $\beta${\rm )}, then $v(0) \in [\xi_0,1]$ {\rm(}or $v(0) \in [-1,\xi_0]${\rm)}.\it
					
					\item If $u(t_0) \in (\alpha, \beta)$, then $v(t)$ is constant in a neighborhood of $t_0$.
					
					\item If $u(t_0) = \alpha$ {\rm(}or $\beta${\rm)}, then $v(t)$ is non-decreasing {\rm(}or non-increasing{\rm)} in a neighborhood of $t_0$.
				\end{enumerate}
		\end{dfn}
	
		By treating $x$ as a parameter, $\hVis (\xi_0, u)$ is defined for $u:Q_T \rightarrow \R$ as we have done previously for $\mathcal{H}(\xi_0,u)$. Visintin\cite{Visintin-1986-evolution} proved the existence of $u$~and~$v$ such that the equation
						\[ u_t = \Delta u + v, \quad v \in \hVis (\xi_0,u), \]
					with $n \geq 1$, Dirichlet boundary conditions, and initial data $\varphi$ is satisfied in a weak sense in $Q_T$.
					Visintin \cite{Visintin-1986-evolution} and more recently Aiki and Kopfova \cite{AikKop-2008-mathematical} proved the existence of solutions to modified versions of \cite{HopJae-1980-pattern,HopJaePoe-1984-hysteresis}, where the hysteretic discontinuity was a completed relay responding to a scalar input. A non-ideal relay with vector input, as in \cite{HopJae-1980-pattern,HopJaePoe-1984-hysteresis}, behaves almost identically to a non-ideal relay with scalar input, but for clarity of exposition we only consider scalar inputs in this chapter.


		Let us now turn to the model hysteresis operator $\hAlt (\xi_0,u)$ proposed by Alt in \cite{Alt-1985-thermostat} (see Figure \ref{fig:hyst-output}c).
		We still consider $H_1(u) \equiv 1$ and $H_{-1} (u) \equiv -1$, and introduce the set
		\[ \tilde{\Sigma_0} := \{ (u,v) \in \R^2 \mid u = \alpha, v \in [-1,1) \} \cup \{ (u,v) \in \R^2 \mid u = \beta, v \in (-1,1] \}\rm{.}\]
			\begin{dfn}\label{dfn-kAlt}
				Let $u,v:[0,T]\rightarrow \R$, where $u$ is a continuous function, and let $\xi_0 \in \{ -1, 1 \}$. We say that  $v \in \hAlt (\xi_0,u)$ if the following hold{\rm:}
					\begin{enumerate}
						\item\label{dfn-Kalt-item1} $(u(t),v(t)) \in \Sigma_{1} \cup \Sigma_{-1} \cup \tilde{\Sigma}_0$ for every $t \in [0,T]$.
						\item\label{dfn-Kalt-item2} If $u(0) \in [\alpha, \beta]$, then $v(0) = \xi_0$.
						\item\label{dfn-Kalt-item3} If $u(t_0) \in (\alpha, \beta)$, then $v(t)$ is constant in a neighborhood of $t_0$.
						\item\label{dfn-Kalt-item4} If $u(t_0) = \alpha$ {\rm(}or $\beta${\rm)}, then $v(t)$ is non-decreasing {\rm(}or non-increasing{\rm)} in a neighborhood of $t_0$.
					\end{enumerate}
			\end{dfn}
One can define $\hAlt (\xi_0,u)$ for $u:Q_T \rightarrow \R$ by treating $x$ as a parameter as we did when defining $\hNir (\xi_0,u)$ and $\hVis (\xi_0,u)$.

		To highlight the main difference between the completed relay\index{completed relay} $\hVis(\xi_0,u)$ and Alt's relay\index{Alt's relay} $\hAlt(\xi_0,u)$, suppose that $\hVis (\xi_0 ,u) (t_0), \hAlt (\xi_0,u)(t_0) \in (-1,1)$ and $u(t_0) = \beta$ has a local maximum at time $t_0$.
		Then, as soon as~$u$ decreases, $\hAlt$ jumps to $-1$, however $\hVis$ remains constant.
		
		Let us introduce the notation $\{ u = \alpha \} := \{ (x,t) \in \overline{Q_T} \mid u(x,t) = \alpha \}$, with $\{ u = \beta \}$ defined analogously. Alt's existence theorem can, omitting the technical assumptions, be stated in the following way. Let $n=1$ and suppose $(\varphi, \xi_0) \in \overline{\Sigma_1} \cup \overline{\Sigma_{-1}}$. Then the following holds:
			\begin{enumerate}
			
				\item  There exists $u$ and $v$ such that $v \in \hAlt(\xi_0,u)$ a.e. in $Q_T$ and
				\[ u_t = u_{xx} + v \quad \text{a.e. on $\{(x,t) \in Q_T \mid u(x,t) \notin \{ \alpha, \beta \} \}$.} \]			
				\item We have		
					\[ u_t = u_{xx} \quad \text{a.e. on $\{ (x,t) \in Q_T \mid u(x,t) \in \{\alpha, \beta \} \}$},\]
				$v \in [-1,0]$ on $\{ u = \beta \}$, and $v \in [0,1]$ on $\{ u = \alpha \}$.
					\item Items \ref{dfn-Kalt-item2}--\ref{dfn-Kalt-item4} of Definition \ref{dfn-kAlt} hold in the following weak sense\index{weak sense}:
					
For every $\psi \in C^\infty_0 (Q \times [0,T))$ with $\psi \geq 0$ on $\{ Q \times [0,T) \} \cap \{ u = \alpha\}$ and $\psi \leq 0$ on $\{ Q \times [0,T) \} \cap \{ u = \beta\}$,
						
						\[ \int_{Q_T} (\eta - \eta_0) \psi_t \, dxdt\leq 0. \]						
					
			\end{enumerate}


\subsection{Slow-Fast Approximation}

\begin{figure}[ht]
\begin{center}
      \includegraphics[width=0.8\linewidth]{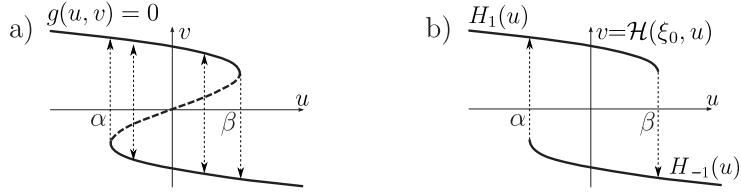}
\end{center}
 \caption{a) The nullcline of the $S$-shaped nonlinearity $g(u,v)$.
b) Hysteresis with nonconstant branches $H_1(u)$ and $H_{-1}(u)$.} \label{fig-SShaped}
\end{figure}

		 Equations of the type \eqref{problem-eq1}--\eqref{problem-eq3} are deeply connected with slow-fast systems\index{slow-fast system} where the variable $v$ is replaced by a fast bistable\index{bistable} ordinary differential equation with a small parameter $\delta > 0$
			 \begin{equation}\label{eqn-slowfast}
			 	\delta v_t=g(u,v).
			 \end{equation} 	 	
		A typical example are the FitzHugh--Nagumo\index{FitzHugh-Nagumo} equations, where $g(u,v) = u + v - v^3 / 3$ and the hysteresis branches\index{hysteresis branches} $H_1 (u)$ and $H_{-1} (u)$ are the stable parts of the nullcline of $g$ (see Figure \ref{fig-SShaped}).
		The question of whether the hysteresis operator approximates the fast variable $v$ as $\delta \rightarrow 0$ has been addressed for systems of ordinary differential equations (see, e.g., \cite{Krejc-2005-hysteresis, MisRoz-1980-differential} and further references in \cite{Kuehn-2015-multiple}), however the corresponding question for partial differential equations is still open.
	\subsection{Free Boundary Approach}		
		Problem \eqref{problem-eq1}--\eqref{problem-eq3} with hysteresis has two distinct phases and a switching mechanism, hence it can be considered as a free boundary problem\index{free boundary problem}. First observe that the hysteresis $\mathcal{H}$ naturally segregates the domain into two subdomains depending on the value of $\xi(x,t)$. Denote			
			\begin{equation}\label{eqn-initconfig}
				Q_j := \{ x \in Q \mid \xi_0(x) = j \}, \quad j= \pm 1.
			\end{equation}			
		Let us look at how the free boundary $\overline{Q_1} \cap \overline{Q_{-1}}$ can evolve for a simple example on the interval $Q = (0,1)$.	
			\begin{figure}[ht!]
			\begin{center}
			\includegraphics[width=1.0\textwidth]{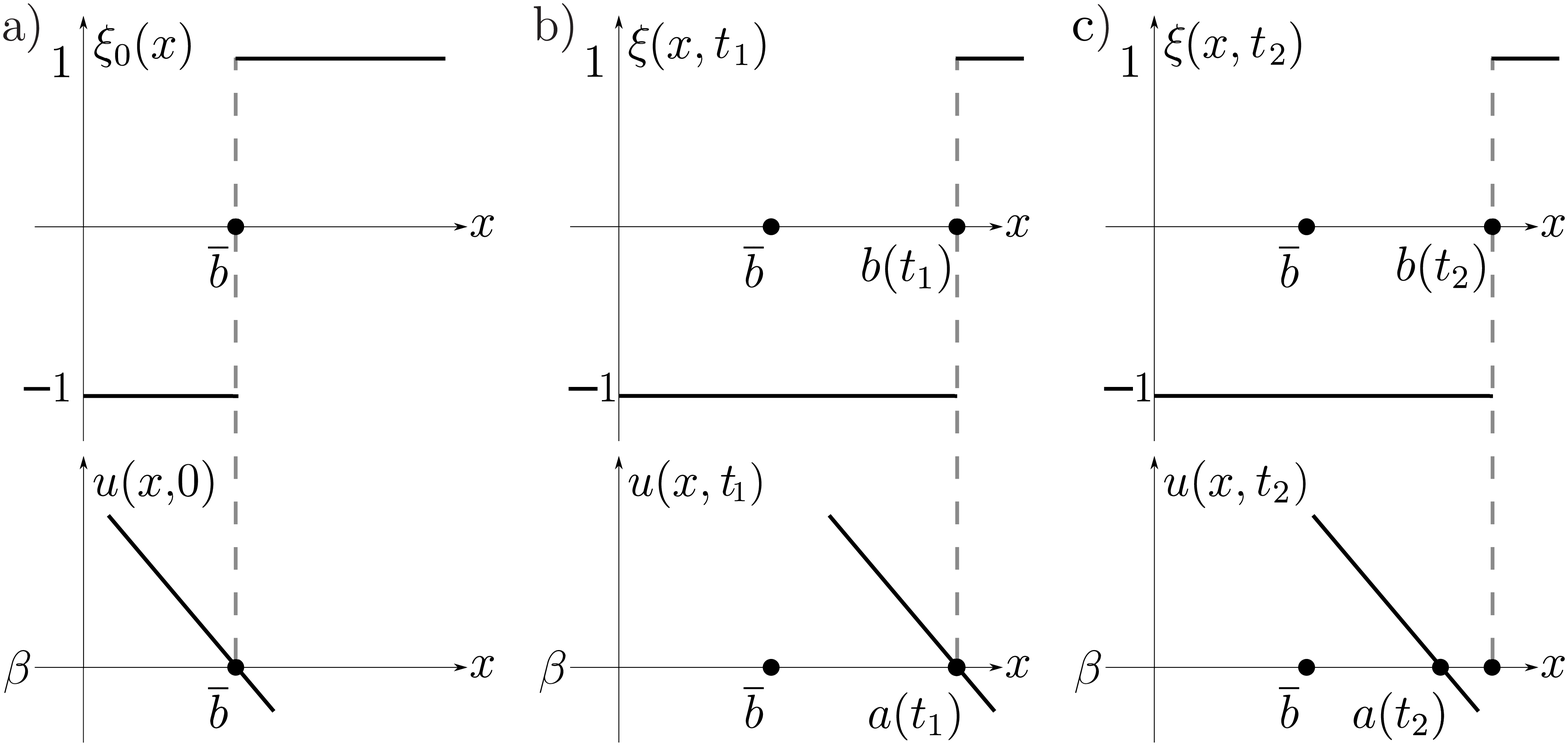}
			\end{center}
			\caption{An example of the hysteresis configuration $\xi$ responding to an input~$u$. }\label{img-root-tracking}
			\end{figure}
		Consider a neighborhood $U$ of $x \in Q$, and suppose at time $t = 0$, $Q_1 \cap U$ and $Q_{-1} \cap U$ are subintervals separated by a point $\overline{b} \in U$ (Figure \ref{img-root-tracking}a).
		Let $u(x,t_0) > \beta$ for $x < \overline{b}$, $u(x,t_0) < \beta$ for $x > \overline{b}$, and let $x = a(t)$ be the unique solution of		
		$u(x,t) = \beta$ in $U$.
		If at time $t_1 > 0$ the value of $u$ at points $x > \overline{b}$ have already risen above $\beta$, then $\xi(x,t)$ has switched from $1$ to $-1$. 	
		These are the points $x$ such that $\overline{b} < x \leq a(t_1)$ (Figure \ref{img-root-tracking}b).
		Now if at time $t_2 > t_1$ the value of $u$ at the switched points has fallen below $\beta$ again, $\xi(x,t)$ remains switched. 	
		These are the points $x$ such that $a(t_2) < x < a(t_1)$ (Figure~\ref{img-root-tracking}c).
		More succinctly, $\xi(x,t) = -1$ if $x \leq b(t)$ and $\xi(x,t) = 1$ if $x > b(t)$, where $b(t) = \max_{0\leq s \leq t} a(s)$.

		The point of this example is to illustrate that the free boundary does not in general coincide with the points where $u$ is equal to one of the threshold values. This is different from the two-phase parabolic obstacle problem (see, e.g., \cite{ApuUra-2014-uniform,ShaUraWei-2009-parabolic}), which \eqref{problem-eq1}--\eqref{problem-eq3} reduces to if $\alpha = \beta$.

		Assume the derivative $\varphi^\prime (\overline{b} )$ in the above example was non-vanishing on the boundary $\{ \overline{b} \} = \overline{Q_1} \cap \overline{Q_{-1}}$.
		This is an example of \emph{transverse}\index{transverse} initial data, and whether the initial data is transverse or not will play an important role in the analysis of problem \eqref{problem-eq1}--\eqref{problem-eq3}.
	
	\subsection{Overview}
		This chapter is organized in the following way.
		
		In Section \ref{sec-transverse} we will investigate the well-posedness of \eqref{problem-eq1}--\eqref{problem-eq3} for \emph{transverse}\index{transverse} initial data. For $n=1$ the existence of solutions and their continuous dependence on initial data was established in \cite{GurTikSha-2013-reaction}, uniqueness of the solution in \cite{GurTik-2012-uniqueness} and the analogous results for systems of equations in \cite{GurTik-2014-systems}.
		Preliminary results for $n \geq 2$ were obtained in \cite{Curran-2014-local}.
		
		In Section \ref{sec-regularity} we consider the regularity\index{regularity} of solutions $u$, in particular, whether the generalized derivatives $u_{x_i x_j}$ and $u_t$ are uniformly bounded. We will summarize the results of \cite{ApuUra-2015-regularity}, where the authors proved that these derivatives are locally bounded in a neighborhood of a point not on the free boundary. They also showed that this bound depends on the parabolic distance to the parts of the free boundary\index{free boundary} that do not contain the sets $\{ u = \alpha \}$ or~$\{ u = \beta \}$.
		

		In Section \ref{sec-non-transverse} we consider non-transverse data and the results of \cite{GurTik-2015-spatially}. We will analyze a spatio-temporal pattern (called \emph{rattling}) arising after spatial discretization of the reaction-diffusion equation and discuss its connection with the continuous model with hysteresis operators $\hNir$, $\hVis$, and $\hAlt$.

	\section{Transverse Initial Data}
	\label{sec-transverse}
		\subsection{Setting of a Model Problem}
			In this section we will discuss the well-posedness\index{well-posedness} of problem \eqref{problem-eq1}--\eqref{problem-eq3} under the assumption that $\varphi$ is transverse\index{transverse} with respect to $\xi_0$, a notion which we will make precise shortly. In order to illustrate the main ideas, we will treat the following model problem in detail and then discuss generalizations at the end of this section (see Subsection \ref{subsec-generalizations}).
			Let $h_{-1} \leq 0 \leq h_1$ be two constants, and let the hysteresis branches be given by $H_1 (u) \equiv h_1$ and $H_{-1} (u) \equiv h_{-1}$.
			Consider the prototype problem
						\begin{align}
							& u_t = \Delta u + \hNir(\xi_0,u), \quad (x,t) \in Q_T,  \label{model-problem-eq1} \\
							& u|_{t=0} = \varphi, \quad x \in Q, \label{model-problem-eq2} \\
							& \frac{\partial u}{\partial \nu} \bigg|_{\partial^\prime Q_T} = 0. \label{model-problem-eq3}
						\end{align}
			We will treat $n =1$ in Subsection \ref{subsec-n=1} (see \cite{GurTikSha-2013-reaction, GurTik-2012-uniqueness}) and $n \geq 2$ (see \cite{Curran-2014-local}) in Subsection \ref{subsec-ngeq2}.		
			Throughout this subsection we will always assume that $\varphi$ and $\xi_0$ are {\it consistent} with each other, i.e., if $\varphi (x) < \alpha$ (or $\varphi(x) > \beta$), then $\xi_0(x) = 1$ (or $\xi_0(x) = -1$).
			In particular, this means that for every $x \in Q$, $\xi(x,t)$ is continuous from the right as a function of $t \in [0,T)$.
  								

	Since in general $\mathcal{H} (\xi_0,u) \in L_q (Q_T)$, we will look for solutions in the Sobolev space $W^{2,1}_q (Q_T)$ with $q > n+2$.	
	This is the space consisting of functions with two weak spatial derivatives and one weak time derivative from $L_q (Q_T)$ (see \cite[Chapter 1]{LadSolUra-1968-Linear}).
	If $ u \in W^{2,1}_q (Q_T)$, then for every $t \in [0,T]$ the trace is well defined and $u(\cdot,t) \in W^{2 - 2 / q}_q (Q)$ (see, e.g., \cite[p. 70]{LadSolUra-1968-Linear}).
To ensure that~$\varphi$ is regular enough to define the spatial transversality property, we henceforth fix a $\gamma$ such that $0 < \gamma < 1 - (n+2) / q$. It follows that if $\varphi \in W^{2 - 2 / q}_q (Q)$, then $\varphi \in C^\gamma(\overline{Q})$ and $\nabla \varphi \in (C^\gamma(\overline{Q}))^n$, where $C^\gamma$ is the standard H\"older space (see \cite[Section 4.6.1]{Triebel-1978-interpolation}).


	The subspace $ W^{2 - 2/q}_{q,N} (Q) \subset  W^{2 - 2 / q}_q (Q)$ of functions with homogeneous Neumann boundary conditions is a well-defined subspace, and in this section we always assume that $\varphi \in W^{2 - 2/q}_{q,N} (Q)$.
	
			\begin{dfn}\label{dfn-soln_existence_problem}
				A solution to problem \eqref{model-problem-eq1}--\eqref{model-problem-eq3} on the time interval $[0,T)$ is a function $u \in
				W^{2,1}_q (Q_T) $ such that
				\eqref{model-problem-eq1} is satisfied in $L_q (Q_T)$ and $u$ satisfies
				\eqref{model-problem-eq2} and \eqref{model-problem-eq3} in terms of traces.
				A solution on $[0,\infty)$ is a function $u:Q \times [0,\infty) \rightarrow \R$ such that for any $T>0$, $u |_{Q_T}$ is a solution in the sense just described.
			\end{dfn}
		We note that if $u \in W^{2,1}_q (Q_T)$, then $\hNir (\xi_0,u)$ is a measurable function on $Q_T$ (see \cite[Section VI.1]{Visintin-1994-differential}).
			
		\subsection{Case \texorpdfstring{$n=1$}{n=1}}
		\label{subsec-n=1}

			Let $Q = (0,1)$ and $Q_j$ be given by \eqref{eqn-initconfig}.
			
				\begin{dfn}\label{dfn-transverse}
					Let $\varphi \in C^1 (\overline{Q})$. We say $\varphi$ is \emph{transverse} with respect to $\xi_0$ if the following hold\emph{:}
						\begin{enumerate}
				\item There is a $\overline{b} \in (0,1)$ such that $ Q_{-1} = \{ x \mid 0 \leq x \leq \overline{b} \}$ and $ Q_1 = \{ x \mid \overline{b} < x \leq 1 \}$.
				\item If $\varphi (\overline{b}) = \beta$, then $\varphi^\prime (\overline{b}) < 0$.
						\end{enumerate}
				\end{dfn}
				An example of $\varphi$ and $\xi_0$ satisfying Definition \ref{dfn-transverse} is given in Figure \ref{img-root-tracking}a.	
					\begin{dfn}\label{dfn-soln-transverse}
						A solution $u$ is called \emph{transverse}\index{transverse} if for all $t \in [0,T]$, $u(\cdot,t)$ is transverse with respect to $\xi(\cdot,t)$.
					\end{dfn}
					
					\begin{thm}[See {\cite[Theorems 2.16 and 2.17]{GurTikSha-2013-reaction}}]\label{thm-local_existence-1d}
						Let $\varphi \in W^{2 - 2/q}_{q,N} (Q)	$ be transverse with respect to $\xi_0$. Then there is a $T > 0$ such that the following hold\emph{:}
							\begin{enumerate}
								\item Any solution $u \in W^{2,1}_q (Q_T)$ of problem \eqref{model-problem-eq1}--\eqref{model-problem-eq3} is transverse.
								\item\label{thm-local-existence-1d-item2} There is at least one transverse solution $u \in W^{2,1}_q (Q_T)$ of problem \eqref{model-problem-eq1}--\eqref{model-problem-eq3}.
								\item If $u \in W^{2,1}_q (Q_T)$ is a transverse solution of problem \eqref{model-problem-eq1}--\eqref{model-problem-eq3}, then it can be continued to a maximal interval of transverse existence $[0,T_{\text{max}})$, i.e., $u(x,T_{\text{max}})$ is not transverse or $T_{\text{max}} = \infty$.
							\end{enumerate}			
					\end{thm}		
			We will sketch the proof of Theorem \ref{thm-local_existence-1d}, part \ref{thm-local-existence-1d-item2}, assuming that $\varphi(\overline{b} ) = \beta$ and $\varphi^\prime (\overline{b}) < 0$.

			Let us define the closed, convex, bounded subset of $C[0,T]$
			\[ B := \{ b \in C[0,T] \mid b(t) \in [0,1], b(0) = \overline{b} \}. \]
			For any $b_0 \in B$, define the function
				\begin{equation}\label{eqn-Fdefinition}	
			 F (x,t) :=
		\left\{ \begin{array}{ll}
		h_{-1} & \text{ if $0 \leq x \leq b_0 (t)$, } \\
		h_1  & \text{ if $b_0(t) < x \leq 1$. }			
		\end{array}  \right.
				\end{equation}
		Let $u \in W^{2,1}_q (Q_T)$ be the solution to problem \eqref{model-problem-eq1}--\eqref{model-problem-eq3} with nonlinearity $F$ in place of $\hNir (\xi_0,u)$.
		We claim that $T$ can be chosen small enough such that the configuration $\xi(x,t)$ of $\mathcal{H}(\xi_0,u)$ is defined by a unique discontinuity point $b(t)$. Note that we do not yet claim that $F = \hNir (\xi_0,u)$.
		
			To prove the claim, first fix $T_0 > 0$. It is a result of classical parabolic theory\index{parabolic theory} \cite[Chapter 4]{LadSolUra-1968-Linear} that for all $T \in [0,T_0]$
				\begin{equation}\label{thm-lin_aux-ineq}
					 \| u \|_{C^\gamma (\overline{Q_T})} + \| u_x \|_{C^\gamma (\overline{Q_T})} \leq C_1 \left( \| F \|_{L_q (Q_T)} + \| \varphi \|_{W_{q,N}^{2 -2 / q}(Q)} \right) \leq C_2,
				\end{equation}	
			where $C_1, C_2, \ldots>0$ depend only on $T_0$ and $q$. The claim now follows from~\eqref{thm-lin_aux-ineq} with the help of the implicit function theorem.
			
			Observe that $u$ is a solution of problem \eqref{model-problem-eq1}--\eqref{model-problem-eq3} if $\hNir (\xi_0,u) = F$, i.e., $b_0 = b$. We therefore look for a fixed point of the map $\mathcal{R}: B \rightarrow B$, $\mathcal{R}(b_0) := b$.
			
			Consider $b_{01}, b_{02} \in B$ and define $F_1, F_2$ via $b_{01}, b_{02}$ similarly to \eqref{eqn-Fdefinition}, and let~$u_1, u_2$ be the corresponding solutions.
			Observe that $F_1 \neq F_2$ only if \linebreak ${{\min (b_{01} (t), b_{02} (t)) < x < \max (b_{01} (t), b_{02} (t))}}$, in particular,
				\begin{equation}\label{eqn-existence-ineq1}		
			\begin{array}{rl}
			\| u_1 - u_2 \|_{C^\gamma (\overline{Q_T})} + \| u_{1x} - u_{2x} \|_{C^\gamma (\overline{Q_T})} & \leq C_1 \| F_{1} - F_{2} \|_{L_q (Q_T)}, \\
				& \leq C_3 \| b_{01} - b_{02} \|_{C[0,T]}^{1 / q}.
			\end{array}						
				\end{equation}
			 Applying \eqref{thm-lin_aux-ineq} again, and using $\varphi^\prime (\overline{b}) > 0$ and the implicit function theorem, we see that
			the left hand side of \eqref{eqn-existence-ineq1} bounds $\| a_1 - a_2 \|_{C[0,T]}$.
			One can additionally show that $\| a_1 - a_2 \|_{C[0,T]}$ bounds $\| b_1 - b_2\|_{C[0,T]}$, hence
				\begin{equation}\label{eqn-1dexistence-1/q-inequality}
			\|b_1 - b_2 \|_{C [0,T]} \leq \|a_1 - a_2 \|_{C [0,T]} \leq C_4 \|b_{01} - b_{02} \|_{C[0,T]}^{1 / q}.
				\end{equation}
			 In particular \eqref{eqn-1dexistence-1/q-inequality} shows that $\mathcal{R}$ is a continuous map on $B$.
			 Moreover, one can use \eqref{thm-lin_aux-ineq} to show that $\mathcal{R}(B)$ is bounded in $C^\gamma [0,T]$, and since $C^\gamma [0,T]$ is compactly embedded into $C[0,T]$, the Schauder fixed point theorem implies that $\mathcal{R}$ has a fixed point.
				\begin{thm}[{see \cite[Theorem 2.2]{GurTik-2012-uniqueness}}]\label{thm-1d-unique}
					If $u_1$ and $u_2$ are transverse solutions of problem \eqref{model-problem-eq1}--\eqref{model-problem-eq3} with the same $\varphi$, then $u_1 \equiv u_2$.
				\end{thm}
				We prove the theorem by expressing solutions as a convolution with the Green function\index{Green function} $G(x,y,t,s)$ for the heat equation with Neumann boundary conditions. Let us use this function to estimate the solution $w= u_1 - u_2$ of the heat equation with zero initial data, Neumann boundary conditions, and the right hand side $h = \mathcal{H}(\xi_0,u_1) - \mathcal{H}(\xi_0,u_2)$:
		\begin{equation}\label{eqn-1d-unique1}
			|w(x,t)| \leq \int_0^t \int_Q |G(x,y,t,s)| | h(y,s)| \, dyds	.	
		\end{equation}
	Also note that $G$ satisfies the inequality (see, e.g., \cite{Ivasishen-1981-Greens})
		\begin{equation}\label{eqn-1d-unique2}
			|G(x,y,t,s)| \leq \frac{C_1}{(t-s)^{1/2}},\quad x,y \in Q, \, 0\leq s < t,		
		\end{equation}	
	where $C_1 > 0$ does not depend on $x$, $y$, $t$ or $s$.
	
	Similarly to the proof of Theorem \ref{thm-local_existence-1d}, for every $s \leq t$ the integral of $|h(y,s)|$ over $Q$ is bounded by $ \| b_1 - b_2 \|_{C[0,t]}$ and hence by $\|a_1 - a_2 \|_{C[0,t]}$ and hence by	$\|u_1 - u_2\|_{C(\overline{Q_t})}$. Combining this with \eqref{eqn-1d-unique1} and \eqref{eqn-1d-unique2}, and taking the supremum over $(x,t) \in Q_T$ we get
		\[ \| w \|_{C(\overline{Q_T})} \leq C_2 \sqrt{T} \| w \|_{C(\overline{Q_T})}, \]
	where $C_2 > 0$ does not depend on $T$. Thus $w = 0$ for $T$ small enough. A passage to arbitrary $T$ is standard.
		\begin{thm}[{See \cite[Theorem 2.9]{GurTikSha-2013-reaction}}]\label{thm-cont-dependence-1d}
			Let $u \in W^{2,1}_q (Q_T)$ be a transverse solution of problem \eqref{model-problem-eq1}--\eqref{model-problem-eq3}.
			If $\| \varphi - \varphi_n \|_{W^{2 - 2/q}_{q,N}(Q)} \rightarrow 0$ and $ | \overline{b}_n - \overline{b}| \rightarrow 0$ as $n \rightarrow \infty$, then for sufficiently large $n$, problem \eqref{model-problem-eq1}--\eqref{model-problem-eq3} has a solution $u_n \in W^{2,1}_q (Q_T)$ with initial data $\varphi_n$ and initial configuration $\xi_{0n}$ defined via~$\overline{b}_n$.
			Furthermore, $\| u_n - u \|_{W^{2,1}_q (Q_T)} \rightarrow 0$ as $n \rightarrow \infty$.
		\end{thm}
		
		The crux of the proof is showing that for sufficiently large $n$, all the solutions exist on the same time interval $[0,T]$.
		To this end we note that we have in fact given an explicit construction of $T$, and that this $T$ depends on~$\overline{b}$, $\| \varphi \|_{W^{2 - 2/q}_q (Q)}$,
		and if $\varphi (\overline{b}) = \beta$, also on $\varphi^\prime (\overline{b})$. Hence for $\varphi_n$ and $\overline{b}_n$ close enough to $\varphi$ and $\overline{b}$ in their respective norms, the same $T$ can be used.

\subsection{Case \texorpdfstring{$n \geq 2$}{neq2}}\label{subsec-ngeq2}

\begin{figure}[ht!]
		\begin{center}
			\includegraphics[width=\textwidth]{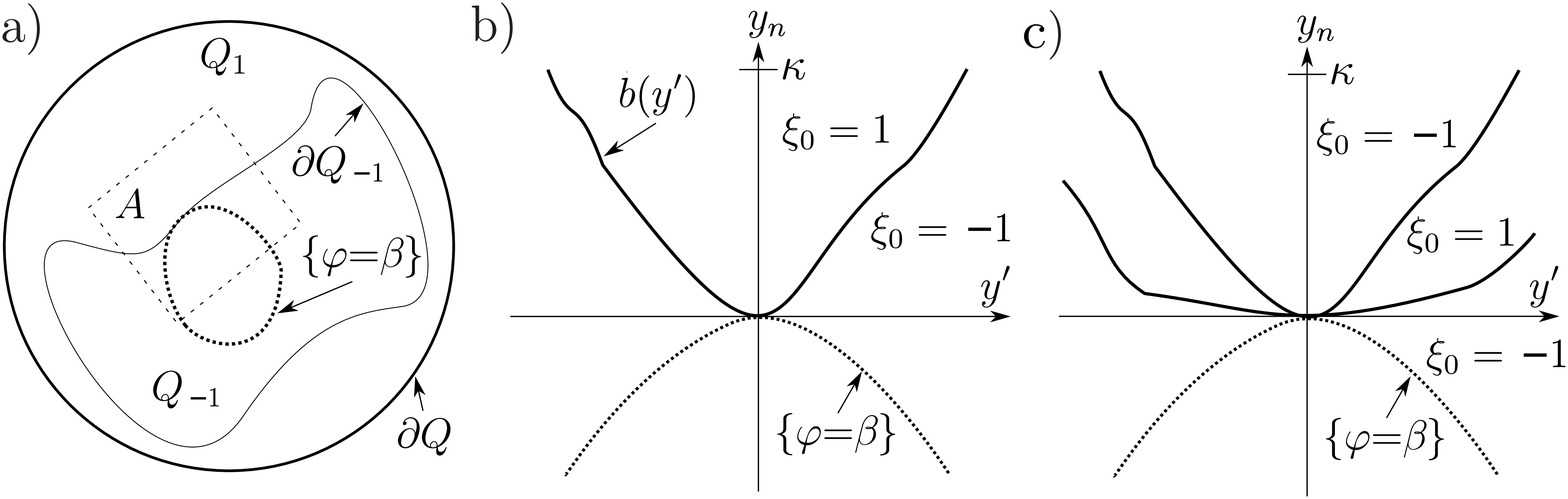}
		\end{center}
			\caption{An example of a) the sets $Q_{\pm1}$, b) transverse data, and c) non-transverse data. \label{img-2d-transverse}}
		\end{figure}

	For the case $ n \geq 2$ a notion of transversality has been studied in a model problem.
	For clarity we will define transversality for the case where the threshold~$\beta$ is adjoined to the free boundary between $ Q_1$ and $Q_{-1}$, and $\alpha$ is not.
	The general case is treated similarly.
	In what follows, let $\text{int}(A)$ denote the topological interior of a subset $A \subset Q$, and let $\{ \varphi = \alpha \}$ be defined similarly to $\{ u = \alpha \}$ but taking $x \in \overline{Q}$ instead of $(x,t) \in \overline{Q_T}$.
	In \cite{Curran-2014-local} the existence and uniqueness of solutions were studied for initial data transverse in the following sense (see Figure \ref{img-2d-transverse}a and \ref{img-2d-transverse}b, and recall that $Q_j$ is given by~\eqref{eqn-initconfig}).
	
		\begin{dfn}\label{dfn-transverse-nd}\index{transverse n $\geq$ 2}
		
					We say the function $\varphi$ is transverse with respect to $\xi_0$ if the following hold{\rm :}
						\begin{enumerate}
							\item\label{dfn-transverse-item1} $Q_1$ and $Q_{-1}$ are measurable, $\partial Q_{-1} \subset Q$, ${\partial Q_1 = \partial Q_{-1} \cup \partial Q}$, and $\partial Q_{-1}$ has zero Lebesgue measure.
							\item\label{dfn-transverse-item2} $\varphi(x) < \beta$ for $x \in \text{\rm int}( Q_1 ) \cup \partial Q$.
							\item\label{dfn-transverse-item3} $\varphi(x) > \alpha$ for $x \in \overline{Q_{-1}}$.
							\item\label{dfn-transverse-graphitem} If $x \in \{ u = \beta \} \cap \partial Q_{-1}$, then there is a neighbourhood $A$ of $x$, a set $A^\prime \subset \R^{n-1}$, a $\kappa >0$, and a map $\psi$ such that
							
								\begin{enumerate}
									\item\label{dfn-transverse-graphitem1} $\psi$ is a composition of a translation and a rotation and
									\[ { \psi(A) = A^\prime \times [-\kappa,\kappa]},\quad \psi(x) = (0,0). \]
									\item\label{dfn-transverse-graphitem2} There is a continuous function $\overline{b}: A^\prime \rightarrow
									 [-\kappa, \kappa]$ such that the configuration function $\xi_0 \circ \psi^{-1}$
									 in  $\psi(A)$ {\rm(}which we denote by $\xi_0(y^\prime,y_n)$, $y^\prime \in A^\prime${\rm)} is given by
									\[ \xi_0(y^\prime,y_n) = \left\{  \begin{array}{ll}
									-1 & \text{ if } -\kappa \leq y_n \leq \overline{b}(y^\prime), \\
									1 & \text{ if } \overline{b}(y^\prime) < y_n \leq \kappa .
									\end{array} \right.
									\]
									\item\label{dfn-transverse-graphitem3} $ \varphi \circ \psi^{-1} $, which we write as $\varphi(y^\prime,y_n)$, satisfies $\varphi_{y_n} (0,0) < 0$.
									\end{enumerate}
						\end{enumerate}
			\end{dfn}
	We observe that in Subsection \ref{subsec-n=1}, the boundary between $Q_1$ and $Q_{-1}$ was a single point $\overline{b}$. But when $n \geq 2$, this boundary is assumed to have the structure of a continuous codimension $1$ submanifold\index{continuous manifold} in a neighborhood of a point on the free boundary\index{free boundary} where $\varphi$ takes a threshold value. Also note that for $n \geq 2$ non-transversality can be caused by the geometry of $\partial Q_{-1}$ in addition to the possible degeneracy of $\nabla \varphi$ (see  Figure \ref{img-2d-transverse}c and Subsection \ref{subsec-generalizations}   for further discussion).
	
		\begin{thm}[see {\cite[Theorems 3.18 and 3.19]{Curran-2014-local}}]\label{thm-nd-wellposedness}
			Assume that $n \geq 2$ and $\varphi \in W^{2 - 2/q}_{q,N} (Q)$ is transverse\index{transverse} with respect to $\xi_0$. Then there is a $T > 0$ such that any solution $u \in W^{2,1}_q (Q_T)$ to problem \eqref{model-problem-eq1}--\eqref{model-problem-eq3} is transverse and there is at least one such solution. Moreover, if for some $T^\prime > 0$, $u_1$ and $u_2$ are two transverse solutions to problem \eqref{model-problem-eq1}--\eqref{model-problem-eq3} on $Q_{T^\prime}$, then $u_1 \equiv u_2$.
		\end{thm}

		The main ideas of the proof are similar to those for the case $n=1$.
		Since $(\varphi(y^\prime,\cdot),\xi_0(y^\prime,\cdot))$ is transverse\index{transverse} in the 1d sense for every $y^\prime \in A^\prime$, one can prove continuity of a map $\mathcal{R}$ that now maps functions $u_0 \in C^\lambda (\overline{Q_T})$ ($\lambda < \gamma$) to solutions $\mathcal{R}(u_0) := u$ of problem \eqref{model-problem-eq1}--\eqref{model-problem-eq3} with the right hand side $\hNir(\xi_0,u_0)$.
		Estimate \eqref{thm-lin_aux-ineq} implies that $u \in C^\gamma (\overline{Q_T})$, and the compactness of the embedding $C^\gamma (\overline{Q_T}) \subset C^\lambda (\overline{Q_T})$ and the Schauder fixed point theorem together imply that $\mathcal{R}$ has a fixed point in $C^\gamma (\overline{Q_T})$.
	

		\subsection{Generalizations and Open Problems}
		\label{subsec-generalizations}

			Let us list some generalizations for the case $n=1$.

				\textbf{Change of topology.} Suppose $u(x,t)$ becomes non-transverse at some time $T$ in the sense of Definition \ref{dfn-transverse}. Then one of two possibilities arise. Either $u(x,T)$ has touched a threshold with zero spatial derivative at some point in $(0,1)$, or this is not the case but $\lim_{t \rightarrow T} b(t) = 1$.
				In the latter case, one can continue the solution, and it remains unique, by redefining the problem effectively without hysteresis \cite[Theorem 2.18]{GurTikSha-2013-reaction}. We say that the topology of the hysteresis\index{hysteresis topology} has changed at time $T$, in the sense that $\xi$ transitions from piecewise constant to uniformly constant.
				
\begin{figure}[H]
				\begin{center}
				\includegraphics[width=0.85\textwidth]{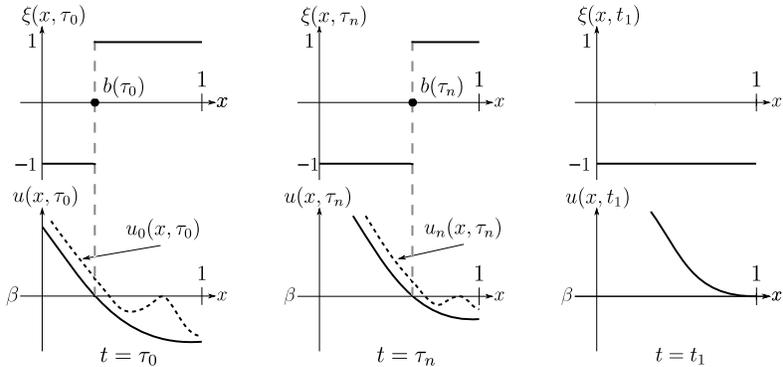}
				\end{center}
				\caption{A solution $u$ (drawn as solid lines in the lower picture) and its configuration $\xi$ (the upper picture) that remain transverse as a discontinuity of $\xi$ disappears at time $t_1$. The dashed line in the lower picture is a series of non-transverse approximations $u_n$ that become non-transverse at moments $\tau_n$ with $\tau_n < t_1$ and $\lim_{n \rightarrow \infty} \tau_n = t_1$. }\label{img-1d-1discontinuity-continuation}
			\end{figure}

\begin{figure}[ht!]
				\begin{center}
				\includegraphics[width=0.85\textwidth]{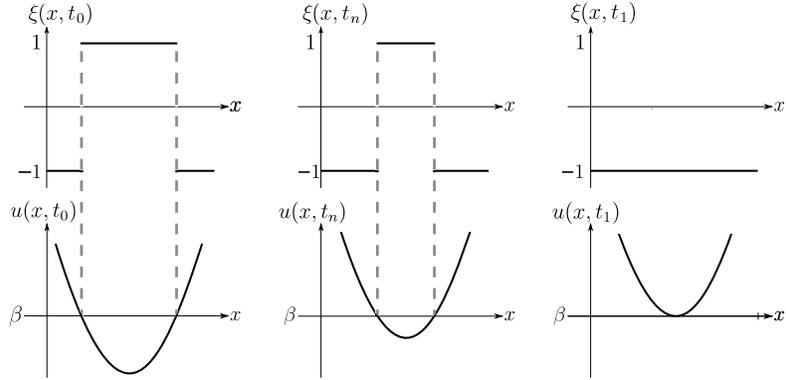}
				\end{center}
				\caption{Discontinuities merging as $t \rightarrow t_1$. }\label{img-1d-discontinuities-merging}
				\end{figure}

				\textbf{Continuous dependence on initial data.}\index{dependence on initial data} If $u$ is a solution such that the topology has changed for some $t_1 < T$, then $u$ need not continuously depend on the initial data since a sequence of approximating solutions $u_n$ may become non-transverse at moments $\tau_n$ with $\tau_n < t_1$ and $\lim_{n \rightarrow \infty} \tau_n = t_1$ (the dashed line in Figure \ref{img-1d-1discontinuity-continuation}).
				But if we also assume that each $u_n$ is a transverse solution, then solutions do depend continuously on their initial data.
				 \textbf{Finite number of discontinuities.} The results in Subsection \ref{subsec-n=1} remain valid if the hysteresis topology is defined by finitely many discontinuity points. The hysteresis changing topology in the sense we described for one point of discontinuity corresponds to these points merging together in the general case (see Figure \ref{img-1d-discontinuities-merging}).
				
				\textbf{General nonlinearity.} The results in this section also hold for the more general problem \eqref{problem-eq1}--\eqref{problem-eq3}. First one must assume that $f$ is locally Lipschitz and dissipative (see \cite[Condition 2.11]{GurTikSha-2013-reaction}). 
				With such an $f$, and if $H_1$ and $H_{-1}$ are locally H\"older continuous, then  transverse solutions exist and can be continuted up to a maximal interval of transverse existence. 
				If one additionally assumes that transverse solutions are unique, they can be shown to continuously depend on their initial data. 
				To prove the uniqueness of solutions the authors of \cite{GurTik-2012-uniqueness,Curran-2014-local} make the stronger assumption on $H_1$ and $H_{-1}$, namely that
							\[ | H_1 (u_1) -H_1 (u_2) | \leq \frac{M}{(\beta - u_1)^{\sigma} + (\beta - u_2)^{\sigma}} |u_1 - u_2 |, \]
							for $u_1, u_2$ in a left neighborhood of $\beta$, with $M > 0$ and $\sigma \in (0,1)$, plus an analogous inequality for $H_{-1}$ and a right neighborhood of $\alpha$. This condition covers the case where $H_{1}$ and $H_{-1}$ are the stable branches in the slow-fast approximation as in Figure \ref{fig-SShaped} (see the appendix of \cite{GurTik-2012-uniqueness} for further discussion).
					
					 \textbf{Systems of equations.} In  {\cite[Theorem 2.1]{GurTik-2014-systems}, the results of Subsection~\ref{subsec-n=1} were generalized to systems of equations of the type in problem \eqref{problem-eq1}--\eqref{problem-eq3}. It was also shown therein that problem \eqref{problem-eq1}--\eqref{problem-eq3} can be coupled to ordinary differential equations to cover the Hoppensteadt--J\"ager model from~\cite{HopJae-1980-pattern,HopJaePoe-1984-hysteresis}.		
		
	Let us conclude this subsection by discussing an open problem.		
		
				\textbf{Open problem.}
	In Figure \ref{img-2d-transverse}c, one can see that for every $y^\prime \neq 0$, \linebreak $(\varphi(y^\prime,\cdot),\xi_0(y^\prime,\cdot))$ is transverse in the 1d sense (with two discontinuties), but since the free boundary cannot be represented as a graph with codomain $y_n$ at the point $y^\prime = 0$, this initial data is not transverse.
				Whether Definition~\ref{dfn-transverse-nd} can be generalized to include such cases is the subject of future work, and at this stage the authors strongly suspect that item \ref{dfn-transverse-graphitem} of Definition \ref{dfn-transverse-nd} can be replaced by the following statement: {\it if $x \in \{u = \beta \} \cap \partial Q_{-1}$, then $\nabla \varphi(x) \neq 0$}.
				In other words, the assumption that the free boundary is a graph is not necessary, and hence Figure~\ref{img-2d-transverse}c would also be transverse.
				This question is intimately linked to the topology of the free boundary.
				Whether solutions can be continued to a maximal interval of existence and how to pose continuous dependence of initial data is unclear for the quite general conditions on $Q_{-1}$~and~$Q_1$ in {Definition \ref{dfn-transverse-nd}}. These questions also apply to the case where $n=1$ and $\xi_0$ has infinitely many discontinuities.

	\section{Regularity of Strong Solutions}\index{regularity of solutions}
	\label{sec-regularity}
	
				To begin with let us discuss what we mean by regularity of solutions in this context. First observe that we cannot expect a classical solution since $\mathcal{H}$ has a jump discontinuity. Therefore the ``optimal'' regularity\index{optimal regularity} we expect is $W^{2,1}_\infty$.
				In this section we obtain $W^{2,1}_\infty$ ``locally'', for points $(x,t) \in Q_T$ outside of the static part of the free boundary\index{static free boundary}. We will also assume the following condition:	
		\begin{con}\label{con-AP-branches}
			$H_1 (u) \equiv 1$ and $H_{-1} (u) \equiv -1$.
		\end{con}	
				 Let us introduce the notation $Q_T^{\pm1} := \{ (x,t) \mid \xi(x,t) = \pm 1 \}$ and observe that~$u$ is smooth on the interior of $Q_T^{\pm1}$.
				
				The free boundary is defined as the set $\Gamma := \partial Q_T^{1} \cap \partial Q_T^{-1}$.
				Moreover, we define $\Gamma_\alpha := \{ u = \alpha \} \cap \Gamma$ and $\Gamma_\beta  := \{ u = \beta \} \cap \Gamma$.
				Note that both $\Gamma_\alpha$ and $\Gamma_\beta$ have zero Lebesgue measure whenever $u$ is a solution of problem \eqref{model-problem-eq1}--\eqref{model-problem-eq3}. This follows from the fact that $u_t - \Delta u = 0$ a.e. on $\Gamma_\alpha \cup \Gamma_\beta$ and Condition \ref{con-AP-branches} (see Alt's argument in the introduction and \cite{Alt-1985-thermostat}).
				
				The estimates we obtain will depend critically on the \emph{static} part of the free boundary $\Gamma_v := \Gamma \backslash (\Gamma_\alpha \cup \Gamma_\beta)$. If $(x,t) \in \Gamma_v$, then $u(x,t) \neq \alpha,\beta$ and by continuity of $u$, $u(x,t \pm \tau) \neq \alpha, \beta$ for $\tau$ sufficiently small. This means $\xi (x,t \pm \tau) = \xi (x,t)$ and so if we draw the $t$-axis vertically as in Figure \ref{fig:Apu}, $\Gamma_v$ looks like a vertical strip.
					\begin{figure}[ht!]
			\centering
			\includegraphics[width=0.5 \textwidth]{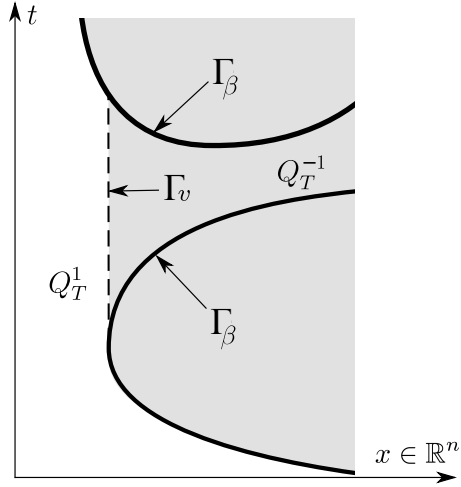}
			\caption{A possible scenario where $\Gamma \neq\Gamma_\alpha \cup \Gamma_\beta$ and $\Gamma_v$ appears. White and grey indicate the regions $Q_T^1$ and $Q_T^{-1}$ respectively.}\label{img-APU-global-picture}\label{fig:Apu}
					\end{figure}	
					
		Next we recall the definition of a parabolic cylinder\index{parabolic cylinder}
		\[ P_r (x^0,t^0) : = \{ x \in \R^n \mid \|x^0 -x \|_{\R^n} < r \} \times (t^0 - r^2, t^0 + r^2), \quad r>0. \]
		We define the parabolic distance between $(x^0,t^0)$ and a set $A \subset Q_T$ as
		\[ \text{dist}_p ((x^0,t^0),A) := \sup \{ r > 0 \mid P_r(x^0,t^0) \cap \{ t \leq t^0\} \cap A  = \emptyset \}. \]
			
			This is all the notation we need to state the main result of \cite{ApuUra-2015-regularity}.
	
			\begin{thm}[see {\cite[Theorem 2.3]{ApuUra-2015-regularity}}]\index{regularity of solutions}
				We assume that $n\geq1$ and $u$ is a solution of problem \eqref{model-problem-eq1}--\eqref{model-problem-eq3}. Then
				\[ |u_t (x,t) | + \sum_{i,j = 1}^n  |u_{x_i x_j} (x,t)| \leq C(\rho_v, \rho_b,M), \quad \text{a.e. ${\rm(}x,t{\rm)} \in Q_T \backslash \overline{\Gamma_v}$}, \]
				where $C$ depends on $\rho_v := \text{dist}_p ((x,t), \Gamma_v )$, $\rho_b := \text{dist}_p ((x,t),  \partial^\prime Q_T \cup ( Q \times \{ 0 \} ))$, and $M:= \sup_{(x,t)\in Q_T}|u(x,t)|$.
			\end{thm}	
%


To explain the main ideas in the proof we define some further notation.
Let $\Gamma_\alpha^0 = \Gamma_\alpha \cap \{ \nabla u = 0\}$ and $\Gamma_\alpha^* = \Gamma_\alpha \backslash \Gamma_\alpha^0$, with $\Gamma_\beta^0$ and $\Gamma_\beta^*$ defined similarly. Furthermore, define $\Gamma^0 = \Gamma_\alpha^0 \cup \Gamma_\beta^0$ and $\Gamma^* = \Gamma_\alpha^* \cup \Gamma_\beta^*$.

The crucial point in the proof is the quadratic growth estimate
\begin{equation} \label{quadratic-growth-u}
\sup\limits_{P_r(x,t)}|u - \beta| \leq C_1(\rho_v,\rho_b,M) r^2 \quad \text{for} \quad r\leq \min \left\lbrace \rho_v, \rho_b \right\rbrace,
\end{equation}
and $(x,t) \in \Gamma^0_\beta$ (the estimate on $\Gamma^0_\alpha$ is similar).
The main tool for showing the  quadratic bound (\ref{quadratic-growth-u}) is the local rescaled version of the Caffarelli monotonicity formula, see \cite{CafSal-2005-geometric,ApuUra-2015-regularity,ApuShaUra-2003-boundary}.

Furthermore, the quadratic growth estimate (\ref{quadratic-growth-u}) implies the corresponding linear bound for $|\nabla u|$
\begin{equation} \label{linear-growth-Du}
\sup\limits_{P_r(x,t)}|\nabla u| \leq C_2(\rho_v,\rho_b,M)r \quad \text{for all}\quad r \leq \min \left\lbrace \rho_v, \rho_b \right\rbrace,
\end{equation}
with $(x,t) \in \Gamma^0$. The dependence of $C_1$ and $C_2$ on the distance $\rho_v$ in \eqref{quadratic-growth-u} and \eqref{linear-growth-Du} arises due to the monotonicity formula. Near $\Gamma_v$ neither  the local rescaled version of Caffarelli's monotonicity formula nor its generalizations (such as the almost monotonicity formula)  are applicable to the positive and negative parts of the spatial directional derivatives $D_e u$, with $e \in \R^n$.

Besides estimates \eqref{quadratic-growth-u} and \eqref{linear-growth-Du}, one also needs  information about the
behaviour of $u_t$ near $\Gamma^*$.
Although $u_t$ may have jumps across the free boundary, one can show that $u_t$ is a continuous function in a neighborhood of $(x,t) \in \Gamma^* \setminus~\Gamma_v$.
In addition, the monotonicity of the jumps of $\hNir(\xi_0,u)$ in the $t$-direction provides  one-sided estimates of $u_t$ near $\Gamma_\alpha$ and $\Gamma_\beta$.
Combining these  results  with the observation that $u_t \leq 0$ on $\Gamma_{\alpha}^* \setminus \Gamma_v$, and $u_t \geq 0$ on $\Gamma_{\beta}^* \setminus \Gamma_v$ gives
\begin{equation} \label{time-derivative-estimate}
\sup_{\Gamma^*\setminus \Gamma_v}|u_t| \leq C_3(\rho_b,M).
\end{equation}
Inequalities (\ref{quadratic-growth-u})--(\ref{time-derivative-estimate}) allow one to apply methods from the theory of free boundary problems (see, e.g., \cite{ApuUra-2014-uniform,ShaUraWei-2009-parabolic}) and  estimate $|u_t(x,t)|$ and $|u_{x_i x_j} (x,t)|$ for a.e. $(x,t) \in Q_T \setminus \overline{\Gamma_v}$.



	\section{Non-Transverse Initial Data}
	\label{sec-non-transverse}
	
	\subsection{Setting of a Problem}
		In this section we summarize the recent work \cite{GurTik-2015-spatially}, where the nontransverse case is analyzed for $x\in\mathbb R$, and indicate directions for further research. We will be interested in the behavior of solutions near one of the thresholds, say~$\beta$. Therefore, we set $\alpha=-\infty$ and $\beta=0$    (see Figure~\ref{figHyst})
\begin{figure}[ht]
\begin{center}
      \includegraphics[width=0.40\linewidth]{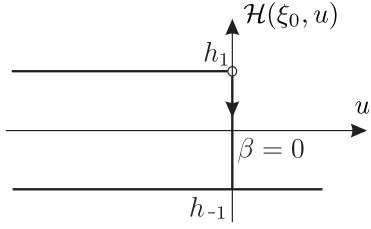}
\end{center}
\caption{Hysteresis with thresholds $\alpha=-\infty$ and $\beta=0$.}
\label{figHyst}
\end{figure}
and assume that the initial data satisfy ${\varphi(x)=-cx^2+o(x^2)}$ in a small neighborhood of the origin, $\varphi(x)<0$ everywhere outside of the origin,   $\xi_0(x)=-1$ for $x=0$, and $\xi_0(x)=1$ for $x\neq 0$. In particular, we assume $c>0$. In this situation, the theorems in Section~\ref{subsec-n=1} are not applicable. Hence, to understand the dynamics of the solution near the origin, we approximate the continuous equation~\eqref{model-problem-eq1} by its spatial discretization and the initial data by the discrete quadratic function. Namely, we choose a grid step $\varepsilon>0$, set $u_n^\varepsilon(t):=u(\varepsilon n,t)$, $n\in\mathbb Z$, and consider the system of infinitely many ordinary differential equations with hysteresis
\begin{equation}\label{eqDiscretePrototype}
\dfrac{d u_n^\varepsilon}{d t}  = \dfrac{u_{n+1}^\varepsilon-2u_n^\varepsilon+u_{n-1}^\varepsilon}{ \varepsilon^2}+\mathcal H(u_n^\varepsilon),\quad t>0, \ n\in\mathbb Z,
\end{equation}
supplemented by the nontransverse (quadratic) initial data
\begin{equation}\label{eqDiscreteInitialData}
  u_n^\varepsilon(0)=  - c(\varepsilon n)^2, \quad n\in\mathbb Z.
\end{equation}
Here we do not explicitly indicate the dependence of $\mathcal H$ on $\xi_0$, assuming that
  $\mathcal H(u_n^\varepsilon)(t)=h_1$ if $u_n^\varepsilon(s)<0$ for all $s\in[0,t]$ and $\mathcal H(u_n^\varepsilon)(t)=h_{-1}$ otherwise. As before, we assume that $h_{-1}\le 0 < h_1$.

Due to~\cite[Theorem~2.5]{GurTik-2015-spatially}, problem~\eqref{eqDiscretePrototype}, \eqref{eqDiscreteInitialData} admits a unique solution in the class of functions satisfying
$$
\sup\limits_{s\in[0,t]}|u_n^\varepsilon(s)|\le A e^{B|n|},\quad n\in\mathbb Z,\ t\ge 0,
$$
with some $A=A(t,\varepsilon)\ge 0$ and $B=B(t,\varepsilon)\in\mathbb R$. Thus, we are now in a position to discuss the dynamics of solutions for each fixed grid step $\varepsilon$ and analyze the limit $\varepsilon\to 0$.

First, we observe that $\varepsilon$ in~\eqref{eqDiscretePrototype}, \eqref{eqDiscreteInitialData} can be scaled out. Indeed, setting
\begin{equation}\label{eqScaleEpsilon}
u_n(t):=\varepsilon^{-2} u_n^\varepsilon(\varepsilon^2 t)
\end{equation}
reduces problem~\eqref{eqDiscretePrototype}, \eqref{eqDiscreteInitialData} to the equivalent one
\begin{equation}\label{eqDiscretePrototypeZEps1}
\left\{
\begin{aligned}
& \dfrac{d u_n}{dt}  = u_{n+1} -2u_n +u_{n-1}  +\mathcal H(  u_n),\quad t>0,\ n\in\mathbb Z,\\
&   u_n(0)=  - c n^2, \quad n\in\mathbb Z.
\end{aligned}\right.
\end{equation}

Using the comparison principle, it is easy to see that if $h_1\le 2c$, then $u_n(t)<0$ for all $n\in\mathbb Z$ and $t>0$ and, therefore, no switchings happen for $t>0$. Let us assume that
\begin{equation}\label{eqh12c}
  h_{-1}\le 0 <2c <h_1.
\end{equation}
It is easy to show that $u_n(t)\le 0$ for all $n\in\mathbb Z$ and $t>0$. However, some nodes can now reach the threshold $\beta=0$ and switch the hysteresis. The main question is which nodes do this and according to which law.

\subsection{Numerical Observations}

 The following pattern formation behavior is indicated by numerics (see Figure~\ref{figRattling}).
\begin{figure}[ht]
\centering
\includegraphics[width=0.95\linewidth]{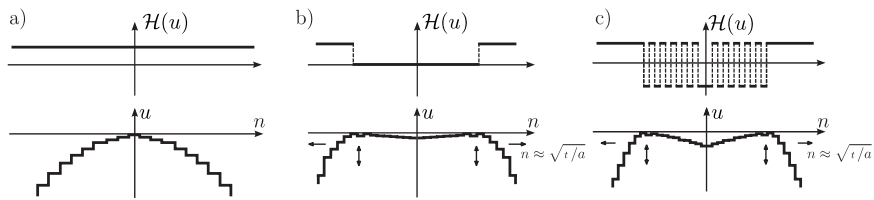}
\caption{Upper graphs represent spatial profiles of the hysteresis
$\mathcal H(u_n)$ and lower graphs the spatial profiles of the solution
$u_n$. a) Nontransverse initial data. b) Spatial profiles at a
moment $t>0$ for $h_{-1}=0$. c) Spatial profiles at a moment
$t>0$ for $h_{-1}=-h_1<0$.}
\label{figRattling}
\end{figure}
As time goes on, the spatial profile of $u_n(t)$ forms two symmetric hills propagating away from the origin.
 At the same time, the whole spatial profile oscillates up and down (never exceeding the threshold $\beta=0$) and touches the threshold $\beta =0$ in such a way that
\begin{equation}\label{eqDiscretePattern}
\lim\limits_{j\to\infty}\dfrac{N_{\rm ns}(j)}{N_{\rm s}(j)}=\dfrac{|h_{-1}|}{h_1},
\end{equation}
where $N_{\rm s}(j)$ and $N_{\rm ns}(j)$ are integers denoting the
number of nodes in the set $\{u_0, u_{\pm1},\dots,u_{\pm j}\}$ that
switch and do not switch, respectively, on the time interval
$[0,\infty)$. In~\cite{GurTik-2015-spatially}, such a spatio-temporal pattern was called {\em rattling}.

A more specific pattern occurs if $|h_{-1}|/h_1=p_{\rm ns}/p_{\rm s}$, where $p_{\rm s}$ and $p_{\rm ns}$ are
co-prime integers. In this case, for any $j$ large enough, the set	
$\{u_{j+1},\dots, u_{j+p_{\rm s}+p_{\rm ns}}\}$ contains exactly $p_{\rm s}$ nodes
that switch and $p_{\rm ns}$ nodes that do not switch on the time
interval $[0,\infty)$.

If a node $u_n$ switches   on the time interval $[0,\infty)$, then we denote its switching moment by  $t_n$; otherwise, set $t_n:=\infty$. In particular, finite values of $t_n$ characterize the propagation velocity of the two hills mentioned above. Numerics indicates that, for the nodes where $t_n$ is finite, we have
\begin{equation}\label{eqtAsymp}
t_n=a n^2+
\begin{cases}
O(\sqrt{  n}) & \text{if } h_{-1}=0,\\
O(  n) & \text{if } h_{-1} < 0,
\end{cases} \quad \text{as } n\to\infty,
 \end{equation}
  and
\begin{equation}\label{eqGradient}
|u_{k+1}(t)-u_{k}(t)|\le b,\quad |k|\le n,\ t\ge t_n,\ n=0,1,2,\dots,
\end{equation}
where $a,b>0$ do not depend on  $k$ and $n$. In particular, \eqref{eqtAsymp} and~\eqref{eqGradient} mean that the hills propagate with velocity of order $t^{-1/2}$, while
the cavity
between the hills has a bounded steepness, which distinguishes the observed phenomenon from the ``classical'' traveling wave situation.

\subsection{Rigorous Result}

The recent work \cite{GurTik-2015-spatially} provides a rigorous analysis of the rattling in the case $h_{-1}=0$, where, according to~\eqref{eqtAsymp}, all the nodes are supposed to switch at time moments satisfying
\begin{equation}\label{eqtnExpected}
t_n=a n^2 + q_n,\quad |q_n|\le E\sqrt{n},
\end{equation}
where $E>0$ does not depend on $n\in\mathbb Z$. In~\cite{GurTik-2015-spatially}, the authors found the coefficient $a$ and proved that if finitely many nodes  $u_n(t)$, $n=0,\pm 1,\dots \pm n_0$, switch at time moments $t_n$ satisfying~\eqref{eqtnExpected}, then all the nodes $u_n(t)$, ${n\in\mathbb Z}$, switch at time moments $t_n$ satisfying~\eqref{eqtnExpected} (see the rigorous statement below). One of the main tools in the analysis is the so-called discrete Green function $y_n(t)$ that is a solution
of the problem
\begin{equation}\label{eqy_nGreen}
\left\{
\begin{aligned}
& \dot y_0=\Delta y_0+1, & & t>0,\\
& \dot y_n=\Delta y_n, & & t>0,\ n\ne0,\\
& y_n(0)=0, & & n\in\mathbb Z.
\end{aligned}
\right.
\end{equation}
The important property of the discrete Green function is the following asymptotics proved in~\cite{Gur-2015-asymptotics}:
\begin{equation}\label{eqAsympyn}
y_n(t)=\sqrt{t} f\!\left(\dfrac{|n|}{\sqrt {t}}\right) + O\left(\dfrac{1}{\sqrt{t}}\right)\quad \text{as } t\to \infty,
\end{equation}
where
\begin{equation}
f(x)  := 2x\int\limits_x^\infty y^{-2}h(y)\,dy,\quad h(x)  := \dfrac{1}{2\sqrt{\pi}}\,e^{-\frac{x^2}{4}},  \quad
\end{equation}
and $O(\cdot)$ does not depend on $n\in\mathbb Z$.

Now if we (inductively) assume that the nodes $u_0, u_{\pm 1},\dots u_{\pm (n-1)}$ switched at the moments satisfying~\eqref{eqtnExpected}, while no other nodes switched on the time interval $[0,t_{n-1}]$, then the dynamics of the node $u_n(t)$ for $t\ge t_{n-1}$ (and until the next switching in the system occurs) is given by
\begin{equation}\label{eqSolu}
u_n(t)=-cn^2+(h_1-2c)t-h_1\sum\limits_{k=-(n-1)}^{n-1} y_{n-k}(t-t_k).
\end{equation}
At the (potential) switching moment $t_n=an^2+q_n$, the relations $t_k=ak^2+q_k$ ($|k|\le n-1$), equality~\eqref{eqSolu}, the Taylor formula, and asymptotics~\eqref{eqAsympyn} yield
\begin{equation}\label{equn1}
\begin{aligned}
0& =  -cn^2+(h_1-2c)an^2-h_1\sum\limits_{k=-(n-1)}^{n-1}  y_{n-k}\left(a(n^2-k^2)\right) + \text{l.o.t.} \\
& = -cn^2+(h_1-2c)an^2-h_1\sum\limits_{k=-(n-1)}^{n-1}  \sqrt{a(n^2-k^2)} f\!\left(\dfrac{n-k}{\sqrt{a(n^2-k^2)}}\right)  \\
& + \text{l.o.t.}\\
& =\left(-c+(h_1-2c)a-h_1 R_n(a)\right)n^2 + \text{l.o.t.},
\end{aligned}
\end{equation}
where
$$
R_n(a):=\sum\limits_{k=-(n-1)}^{n-1}\dfrac{1}{n}  \sqrt{a(1-(k/n)^2)} f\!\left(\dfrac{1-k/n}{\sqrt{a(1-(k/n)^2)}}\right)
$$
and ``l.o.t.'' stands for lower order terms that we do not explicitly specify here. Note that $R_n(a)$ is the Riemann sum for the integral
\begin{equation}\label{eqIf}
I_f(a):=\int\limits_{-1}^1  \sqrt{a(1-x^2)} f\!\left(\dfrac{1-x}{\sqrt{a(1-x^2)}}\right)\, dx.
\end{equation}
Therefore, equality~\eqref{equn1} can be rewritten as
\begin{equation}\label{equn2}
0 = \left(-c+(h_1-2c)a-h_1 I_f(a)\right)n^2 + \text{l.o.t.}
\end{equation}
It is proved in~\cite{GurTik-2015-spatially} that there exists a unique $a>0$ for which the coefficient at $n^2$ in~\eqref{equn2} vanishes. The most difficult part is to analyze the lower order terms in~\eqref{equn2} that involve:
\begin{enumerate}
\item the remainders $q_0,q_{\pm1},\dots,q_n$ from~\eqref{eqtnExpected} arising from~\eqref{eqSolu} via the application of the Taylor formula,

\item the remainder in the asymptotic~\eqref{eqAsympyn} for  the discrete Green function $y_n(t)$,

\item the remainders arising from approximating the integral $I_f(a)$ by the Riemann sum $R_n(a)$.
\end{enumerate}
In particular, one has to prove that if $|q_j|\le E\sqrt{|j|}$ for $j=0,\pm1,\dots,\pm (n-1)$, then the lower order terms vanish for $a$ specified above and $|q_n|\le E\sqrt{|n|}$. This allows one to continue the inductive scheme and (after an appropriate analysis of the nodes $u_{\pm(n+1)}(t), u_{\pm(n+2)}(t),\dots$ for $t\in[t_{n-1},t_n]$) complete the proof.

The rigourous formulation of the main result in \cite{GurTik-2015-spatially} is as follows.
\begin{thm}[see~{\cite[Theorem~3.2]{GurTik-2015-spatially}}]
Assume that \eqref{eqh12c} holds and that \linebreak ${h_{-1}=0}$. Let $a=a(h_1/c)>0$ be a $($unique$)$ root of the equation
\begin{equation}\label{eqRoots}
-c+(h_1-2c)a-h_1 I_f(a)=0
\end{equation}
with $I_f(a)$ given by~\eqref{eqIf}. Then there exists a constant $E_0=E_0(h_1,c,a)>0$ and a function $n_0=n_0(E)=n_0(E,h_1,c,a)$ $($both explicitly constructed$)$ with the following property. If
\begin{equation}\label{eqFinitelyManySwitchings}
\begin{aligned}
 & \text{finitely many nodes}\ u_0(t),u_1(t),\dots,u_{n_0}(t)\ \text{switch at  moments}\ t_n\\
 &
\text{satisfying~\eqref{eqtnExpected} with the above $a$ and some $E\ge E_0$},
\end{aligned}
\end{equation}
then each node $u_n(t)$, $n\in\mathbb Z$,
switches{\rm;} moreover, the switching occurs at a time moment $t_n$
satisfying~\eqref{eqtnExpected} with $a$ and $E$ as in~\eqref{eqFinitelyManySwitchings}.
\end{thm}

We note that the explicit formula~\eqref{eqSolu} for the solution $u_n(t)$ allows one to verify
the fulfillment of finitely many
assumptions~\eqref{eqFinitelyManySwitchings}
numerically with an arbitrary accuracy for any given values of $h_1$
and $c$. The graphs in Figure~\ref{figaEn0} taken from~\cite{GurTik-2015-spatially} represent the values of $a$, $E$, and
$n_0(E)$ that fulfill assumption~\eqref{eqFinitelyManySwitchings}
for $c=1/2$ and $h_1=1.1, 1.2, 1.3,  \dots, 2.5$.
\begin{figure}[ht]
\begin{center}
      \includegraphics[width=0.92\linewidth]{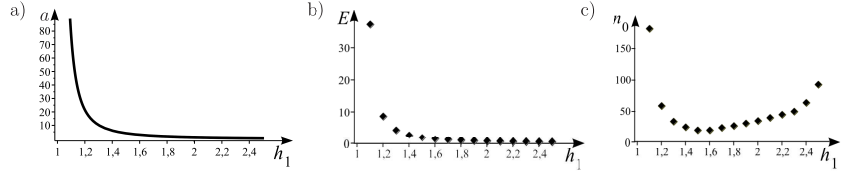}
\end{center}
\caption{Dependence on $h_1$ of the values of $a$, $E$, and
$n_0(E)$ that fulfill assumptions~\eqref{eqFinitelyManySwitchings}
for $c=1/2$. a) The values of $a$ are found as roots of~\eqref{eqRoots}. b),~c)~The values of $E$ and
$n_0(E)$ are calculated for discrete values $h_1=1.1, 1.2, 1.3,  \dots, 2.5$.}
\label{figaEn0}
\end{figure}

\subsection{Open Problems}

To conclude this section, we indicate several directions of further research in the nontransverse case.

{\bf Case $h_{-1}<0$.} In this case, one has to additionally prove a specific switching pattern~\eqref{eqDiscretePattern}. We expect that the tools developed in~\cite{GurTik-2015-spatially} will work for rational $h_1/h_{-1}$. The irrational case appears to be a much more difficult problem.

{\bf Multi-dimensional case.} Numerics indicates that the behavior analogous to~\eqref{eqDiscretePattern} occurs in higher spatial dimensions for different kinds of approximating grids. Figure~\ref{figMultiDim} illustrates the switching pattern
for a two-dimensional analog of problem~\eqref{eqDiscretePrototypeZEps1}, where the Laplacian is discretized on the square
and triangular lattices, respectively.
\begin{figure}[ht]
\begin{center}
      \includegraphics[width=0.70\linewidth]{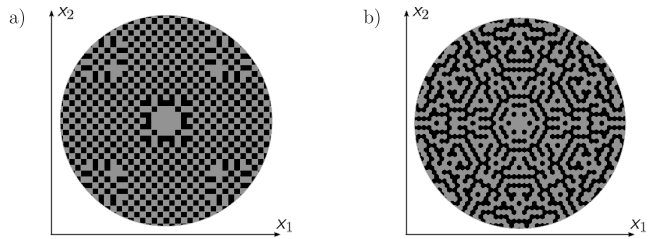}
\end{center}
\caption{A snapshot for a time moment $t>0$ of a
two-dimensional spatial profile of hysteresis taking values
$h_1>2c>0$ and $h_{-1}=-h_1<0$. The nontrasverse initial data is
given by $\varphi(x)=-c(x_1^2+x_2^2)$. Grey (black) squares or
hexagons correspond to the nodes that have (not) switched on the
time interval $[0,t]$. a) Discretization on the square lattice.
b) Discretization on the triangular lattice.} \label{figMultiDim}
\end{figure}

{\bf Limit $\varepsilon\to 0$}. We introduce the function
$$
u^\varepsilon(x,t):=u_n^\varepsilon(t),\quad
 x\in[\varepsilon n-\varepsilon/2,\varepsilon n + \varepsilon/2),\ n\in\mathbb Z,
$$
(which is piecewise constant in $x$ for every fixed $t$). Making the transformation inverse to~\eqref{eqScaleEpsilon} and assuming  \eqref{eqDiscretePattern} and \eqref{eqtAsymp}, we can deduce that, as $\varepsilon\to 0$, the function $u^\varepsilon(x,t)$ approximates a smooth function $u(x,t)$, which satisfies
$u(x,t)=0$ for $x\in(-\sqrt{t/a},\sqrt{t/a})$. In other words, $u(x,t)$ sticks to the
threshold line $\beta=0$ on the expanding interval
$x\in(-\sqrt{t/a},\sqrt{t/a})$.

Similarly to $u^\varepsilon(x,t)$, we consider the function
$$
H^\varepsilon(x,t):=\mathcal H(u_n^\varepsilon)(t),\quad
x\in[\varepsilon n-\varepsilon/2,\varepsilon n + \varepsilon/2), \
n\in\mathbb Z,
$$
which is supposed to approximate the hysteresis
$\mathcal H(u)(x,t)$ in~\eqref{model-problem-eq1}.
 We see that the spatial profile of $H^\varepsilon(x,t)$ for $x\in(-\sqrt{t/a},\sqrt{t/a})$ is a step-like function taking values $h_1$ and $h_{-1}$ on alternating intervals of length of order $\varepsilon$. Hence, it
has no pointwise limit as $\varepsilon\to 0$, but converges
in a weak sense to the function $H(x,t)$ given by
$H(x,t)=0$ for $x\in(-\sqrt{t/a},\sqrt{t/a})$ and $H(x,t)=h_1$ for $x\notin(-\sqrt{t/a},\sqrt{t/a})$.
We emphasize that $H(x,t)$ does not
depend on $h_{-1}$ (because $a$ does not). On the other hand, if
$h_{-1} < 0$, the hysteresis operator $\mathcal H(u)(x,t)$
in~\eqref{model-problem-eq1} cannot take value $0$ by definition, which
clarifies the essential difficulty with the well-posedness of the
original problem~\eqref{model-problem-eq1}  in the nontransverse case. To
overcome the non-wellposedness, one need to allow the intermediate
value $0$ for the hysteresis operator, cf. the discussion of modified hysteresis operators due to Visintin and Alt in the introduction. A rigorous analysis of the limit  $\varepsilon\to 0$ is an open problem, which may lead to a unique ``physical'' choice of an appropriate element in the multi-valued Visintin's hysteresis $\mathcal H_{\rm Vis}(\xi_0,u)$ in Definition~\ref{dfn-kVis}.

{\bf Rattling in slow-fast systems.} One may think that the rattling occurs exclusively due to the discontinuous nature or hysteresis. This is not quite the case. Consider an equation of type~\eqref{model-problem-eq1} with the hysteresis $\mathcal H(\xi_0,u)$ replaced by the solution $v$ of a bistable ordinary differential equation of type~\eqref{eqn-slowfast}, e.g.,
\begin{equation}\label{eqn-completeSF}
u_t=u_{xx}+v,\quad \delta v_t=g(u,v).
\end{equation}

Numerical solution of system~\eqref{eqn-completeSF} with a nontransverse initial data $u(x,0)=-cx^2+o(x^2)$ and $v(x,0)=H_1(\beta)$ near the origin reveals a behavior analogous to that for a spatially discrete system (see Figure~\ref{figSlowFastPattern}). As the spatial profile of $u(x,t)$ touches the threshold $\beta$ at some point $x_0$, the spatial profile of $v(x,t)$ forms a peak-like transition layer around $x_0$ that rapidly converges to a plateau. Thus, as time goes on, the spatial profile of $v(x,t)$  converges to a step-like function  taking values $H_{1}(\beta)$ and $H_{-1}(\beta)$ on alternating intervals, whose length tends to zero as $\delta\to 0$. A rigorous analysis of the limit  $\delta\to 0$ is an open problem.
\begin{figure}[ht]
\begin{center}
      \includegraphics[width=0.70\linewidth]{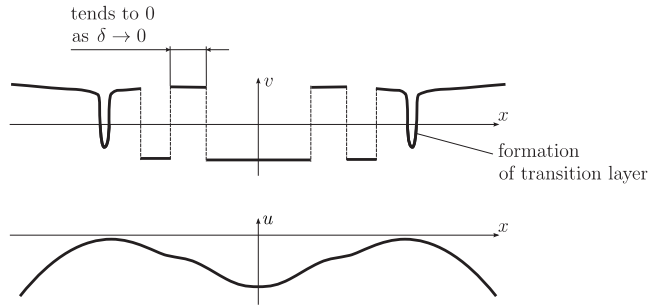}
\end{center}
\caption{Lower and upper graphs are spatial profiles of the solution
$u(x,t)$ and $v(x,t)$, respectively, for
problem~\eqref{eqn-completeSF} with initial data
$u|_{t=0}=-cx^2 + o(x^2)$, $v|_{t=0}=H_1(\beta)$.} \label{figSlowFastPattern}
\end{figure}

\subsubsection*{Acknowledgements}

The authors are grateful for the support of the DFG project SFB 910 and the DAAD project G-RISC. The work of the first author was partially supported by the Berlin Mathematical School. The work of the second author was partially supported by the DFG Heisenberg Programme. The work of the third author was partially supported by Chebyshev Laboratory (Department of Mathematics and Mechanics, St. Petersburg State University)  under RF Government grant 11.G34.31.0026, JSC ``Gazprom neft'', by the Saint-Petersburg State University research grant 6.38.223.2014 and RFBR 15-01-03797a.

\bibliographystyle{ieeetr}
\bibliography{spatially-distributed-hysteresis-spphys}		

\begin{thebibliography}{10}

\bibitem{Visintin-1994-differential}
A.~Visintin, {\em Differential Models of Hysteresis}.
\newblock Applied Mathematical Sciences, Berlin Heidelberg: Springer-Verglag,
  1994.

\bibitem{krejc-1996-hysteresis}
P.~Krej{\v{c}}{\'\i}, {\em Hysteresis, Convexity and Dissipation in Hyperbolic
  Equations}.
\newblock GAKUTO International series, Gatt{\"o}toscho, 1996.

\bibitem{BroSpr-1996-hysteresis}
M.~Brokate and J.~Sprekels, {\em Hysteresis and Phase Transitions}.
\newblock Applied Mathematical Sciences, New York: Springer-Verlag, 1996.

\bibitem{Visintin-2014-ten}
A.~Visintin, ``Ten issues about hysteresis,'' {\em Acta Applicandae
  Mathematicae}, vol.~132, no.~1, pp.~635--647, 2014.

\bibitem{Visintin-2015-P.D.E.s}
A.~Visintin, ``{P.D.E.s with hysteresis 30 years later.},'' {\em Discrete
  Contin. Dyn. Syst., Ser. S}, vol.~8, no.~4, pp.~793--816, 2015.

\bibitem{HopJae-1980-pattern}
F.~Hoppensteadt and W.~J{\"a}ger, ``{Pattern Formation by Bacteria},'' in {\em
  {Biological Growth and Spread}} (W.~J{\"a}ger, H.~Rost, and P.~Tautu, eds.),
  vol.~38 of {\em {Lecture Notes in Biomathematics}}, pp.~68--81, Springer
  Berlin Heidelberg, 1980.

\bibitem{HopJaePoe-1984-hysteresis}
F.~Hoppensteadt, W.~J{\"a}ger, and C.~P{\"o}ppe, ``{A hysteresis model for
  bacterial growth patterns},'' in {\em {Modelling of Patterns in Space and
  Time}} (W.~J{\"a}ger and J.~D. Murray, eds.), vol.~55 of {\em {Lecture Notes
  in Biomathematics}}, pp.~123--134, Springer Berlin Heidelberg, 1984.

\bibitem{Marciniak-2006-receptor}
A.~Marciniak-Czochra, ``Receptor-based models with hysteresis for pattern
  formation in hydra,'' {\em Mathematical Biosciences}, vol.~199, no.~1, pp.~97
  -- 119, 2006.

\bibitem{Koethe-2013-hysteresis}
A.~K{\"o}the, {\em Hysteresis-Driven Pattern Formation in
  Reaction-Diffusion-ODE Models}.
\newblock PhD thesis, University of Heidelberg, 2013.

\bibitem{Krasnosel-2012-systems}
M.~Krasnosel'skii, M.~Niezgodka, and A.~Pokrovskii, {\em Systems with
  Hysteresis}.
\newblock Springer Berlin Heidelberg, 2012.

\bibitem{GurTikSha-2013-reaction}
P.~Gurevich, S.~Tikhomirov, and R.~Shamin, ``Reaction diffusion equations with
  spatially distributed hysteresis,'' {\em Siam J. of Math. Anal.}, vol.~45,
  no.~3, pp.~1328--1355, 2013.

\bibitem{Alt-1985-thermostat}
H.~W. Alt, ``On the thermostat problem,'' {\em Control Cybern.}, vol.~14,
  no.~1-3, pp.~171--193, 1985.

\bibitem{Visintin-1986-evolution}
A.~Visintin, ``Evolution problems with hysteresis in the source term,'' {\em
  SIAM J. Math. Anal.}, vol.~17, no.~5, 1986.

\bibitem{AikKop-2008-mathematical}
T.~Aiki and J.~Kopfov{\'a}, ``A mathematical model for bacterial growth
  described by a hysteresis operator,'' in {\em Recent Advances in Nonlinear
  Analysis}, pp.~1--10, 2008.

\bibitem{Krejc-2005-hysteresis}
P.~Krej\v{c}\'i, ``The hysteresis limit in relaxtion oscillation problems,''
  {\em J. Physics.: Conf. Ser.}, no.~22, pp.~103--123, 2005.

\bibitem{MisRoz-1980-differential}
E.~Mischenko and N.~Rozov, {\em Differential Equations with Small Parameters
  and Relaxation Oscillations}.
\newblock New York: Plenum, 1980.

\bibitem{Kuehn-2015-multiple}
C.~Kuehn, {\em Multiple Time Scale Dynamics}, vol.~191 of {\em Applied
  Mathematical Sciences}.
\newblock Springer International Publishing, 2015.

\bibitem{ApuUra-2014-uniform}
D.~Apushkinskaya and N.~Uraltseva, ``{Uniform estimates near the initial state
  for solutions of the two-phase parabolic problem.},'' {\em St. Petersbg.
  Math. J.}, vol.~25, no.~2, pp.~195--203, 2014.

\bibitem{ShaUraWei-2009-parabolic}
H.~Shahgholian, N.~Uraltseva, and G.~S. Weiss, ``A parabolic two-phase
  obstacle-like equation,'' {\em Adv. Math.}, vol.~221, no.~3, pp.~861--881,
  2009.

\bibitem{GurTik-2012-uniqueness}
P.~Gurevich and S.~Tikhomirov, ``Uniqueness of transverse solutions for
  reaction-diffusion equations with spatially distributed hysteresis,'' {\em
  Nonlinear Anal.}, vol.~75, pp.~6610--6619, December 2012.

\bibitem{GurTik-2014-systems}
P.~Gurevich and S.~Tikhomirov, ``Systems of reaction-diffusion equations with
  spatially distributed hysteresis.,'' {\em Mathematica Bohemica (Proc.
  Equadiff 2013)}, vol.~139, no.~2, pp.~239--257, 2014.

\bibitem{Curran-2014-local}
M.~Curran, ``Local well-poseness of a reaction-diffusion equation with
  hysteresis.,'' Master's thesis, Fachbereich Mathematik und Informatik, Freie
  Universit{\"a}t Berlin, 2014.

\bibitem{ApuUra-2015-regularity}
D.~Apushkinskaya and N.~Uraltseva, ``On regularity properties of solutions to
  the hysteresis-type problem,'' {\em Interfaces and Free Boundaries}, vol.~17,
  no.~1, pp.~93--115, 2015.

\bibitem{GurTik-2015-spatially}
P.~{Gurevich} and S.~{Tikhomirov}, ``Spatially discrete reaction-diffusion
  equations with discontinuous hysteresis,'' {\em ArXiv:1504.02385 [math.AP]},
  2015.

\bibitem{LadSolUra-1968-Linear}
O.~Ladyzhenskaya, V.~Solonnikov, and N.~Uraltseva, {\em Linear and Quasilinear
  Equations of Parabolic Type}.
\newblock Providence, Rohde Island: American Mathematical Society, 1968.

\bibitem{Triebel-1978-interpolation}
H.~Triebel, {\em Interpolation Theory, Function Spaces, Differential
  Operators}.
\newblock Carnegie-Rochester Conference Series on Public Policy, North-Holland
  Publishing Company, 1978.

\bibitem{Ivasishen-1981-Greens}
S.~Ivasishen, ``Green's matrices of boundary value problems for {P}etrovskii
  parabolic systems of general form. ii,'' {\em Math. USSR-Sb}, no.~4,
  pp.~461--498, 1981.

\bibitem{CafSal-2005-geometric}
L.~Caffarelli and S.~Salsa, {\em A Geometric Approach to Free Boundary
  Problems}.
\newblock Graduate studies in mathematics, American Mathematical Soc., 2005.

\bibitem{ApuShaUra-2003-boundary}
D.~Apushkinskaya, H.~Shahgholian, and N.~Uraltseva, ``Boundary estimates for
  solutions of the parabolic free boundary problem,'' {\em Journal of
  Mathematical Sciences}, vol.~115, no.~6, pp.~2720--2730, 2003.

\bibitem{Gur-2015-asymptotics}
P.~{Gurevich}, ``{Asymptotics of parabolic Green's functions on lattices},''
  {\em ArXiv:1504.02673 [math.AP]}, 2015.

\end{thebibliography}

\end{document}